\documentclass[11pt, reqno]{amsart}
\usepackage[utf8]{inputenc}
\usepackage{amsmath}
\usepackage{amsthm}
\usepackage{amsfonts}
\usepackage{amssymb}
\usepackage{mathrsfs, mathtools, mathdots}
\usepackage{diagbox}
\usepackage{relsize}

\usepackage{enumerate, xcolor, graphicx}
\definecolor{Magenta}{cmyk}{0,1,0,0}
\definecolor{dgreen}{rgb}{0.,0.6,0.}

\usepackage[top=27mm,bottom=27mm,left=25mm,right=25mm]{geometry}

\newtheorem{theorem}{Theorem}[section]

\newtheorem{corollary}{Corollary}[section]

\newtheorem*{definition*}{Definition}
\newtheorem*{theorem*}{Theorem}
\newtheorem*{remark*}{Remark}
\newtheorem*{problem*}{Problem}
\newtheorem*{conjecture*}{Conjecture}

\newcommand{\Q}{\mathbb{Q}}
\newcommand{\nq}{[n]_{q}}

\begin{document}

\title[Values of the $q$-exponential and related functions]{Linear independence of values of \\ the $q$-exponential and related functions}

\author[A. B. Dixit]{Anup B. Dixit}

\author[V. Kumar]{Veekesh Kumar}

\author[S. Pathak]{Siddhi S. Pathak}

\address{Institute of Mathematical Sciences (HBNI),  CIT Campus Taramani, Chennai, Tamil Nadu, India 600113}

\address{Department of Mathematics, Indian Institute of Technology, Dharwad, Karnataka, India 580011.}

\address{Chennai Mathematical Institute, H-1 SIPCOT IT Park, Siruseri, Kelambakkam, Tamil Nadu, India 603103.}

\email{anupdixit@imsc.res.in}

\email{veekeshk@iitdh.ac.in}

\email{siddhi@cmi.ac.in}

\date{\today}



\thanks{Research of the first and the third authors was partially supported by the INSPIRE Faculty Fellowship.}

\begin{abstract}
In this paper, we establish the linear independence of values of the $q$-analogue of the exponential function, $E_q(x)$ and its derivatives at specified algebraic arguments, when $q$ is a Pisot-Vijayraghavan number. We also deduce similar results for cognate functions, such as the Tschakaloff function and certain generalized $q$-series.

\end{abstract}

\maketitle

\section{\bf Introduction}
\bigskip

For any complex number $q$ with $|q| > 1$, the $q$-analogue of the exponential function is defined by the absolutely convergent series
\begin{equation*}
    E_q(x) := 1+ \sum_{n=1}^{\infty} \frac{x^n}{\nq !},
\end{equation*}
where $\nq = (q^n-1)$ and $\nq! = (q^n-1)(q^{n-1}-1) \cdots (q-1)$. Similarly we have the $q$-analogue of the logarithm, which is given by 
\begin{equation*}
    L_q(x) := \sum_{n=1}^{\infty} \frac{x^n}{\nq}, \qquad \text{for } |x| < |q|.
\end{equation*}
The analogy between the classical functions and their $q$-analogues is driven by the limit
\begin{equation*}
    \lim_{q \rightarrow 1^+} \frac{q^n - 1}{q-1} = n.
\end{equation*}
Unlike the classical exponential and logarithm functions, their $q$-counterparts are related by the following differential relation:
\begin{equation*}
    L_q(x) = x \frac{E'_q(-x)}{E_q(-x)}
\end{equation*}
for $|x| < |q|$. For more details, we refer the reader to \cite[Section 6]{murty-fib}. These functions appear in various contexts in combinatorics and number theory, and are studied as interesting functions in their own right. \\

The value at $x=1$ of the $q$-logarithm function is of particular importance, as $L_q(1) = \zeta_q(1)$, where
\begin{equation*}
    \zeta_q(s) := \sum_{n=1}^{\infty} \frac{n^{s-1}}{\nq},
\end{equation*}
is the $q$-analogue of the Riemann zeta-function as considered in \cite{krz}. The value $\zeta_q(1)$,
\begin{equation*}
    \zeta_q(1) = \sum_{n=1}^{\infty} \frac{1}{q^n -1},
 \end{equation*}
is often referred to in the literature as the $q$-harmonic series.\\

In this paper, we examine the arithmetic nature and linear independence properties of certain special values of the above mentioned functions. Recall that a real algebraic integer $\omega$ is said to be \textit{Pisot-Vijayraghavan number} (abbreviated as PV number) if $\omega > 1$, and for all other Galois conjugates $\omega^{(j)}$ of $\omega$, we have $|\omega^{(j)}| < 1$.  Immediate examples of PV numbers are positive integers greater than $1$. A non-trivial example is obtained by considering $\beta$, the real root of $x^4 -x^3 -2x^2 +1$ with $\beta > 1$. Then it can be checked that $\beta$ is a PV number. In fact, Pisot \cite{pisot} showed that in every real algebraic number field, there exist PV numbers that generate the field. These numbers make a fundamental appearance in Diophantine approximation and have been extensively studied in the literature. \\

Fix an algebraic integer $q \neq 0$ and let $n_q = [\Q(q) : \Q]$. Let $ \sigma_1$, $\sigma_2$, $\cdots$, $\sigma_{n_q}$ denote the embeddings of $\Q(q)$ into $\mathbb{C}$, with $\sigma_1$ being identity. Let $\mathcal{O}_{q}$ be the ring of integers of $\Q(q)$. For any algebraic number $\alpha \in \Q(q)$, the $q$-relative height of $\alpha$, $H_q(\alpha)$, is defined as
\begin{equation*}
    H_q(\alpha) := \prod_{l=1}^{n_q} \max  \left\{ 1, \,  \left| \sigma_l(\alpha) \right| \right\}.
\end{equation*}
Thus, if $q$ is a PV number, then $H_q(q) = q$.\\

Our first theorem concerns the linear independence of values of derivatives of a certain generalized $q$-exponential function. Let $P(X) \in \mathbb{Z}[X]$ be a non-constant polynomial such that $P(q^t) \neq 0$ for all $t \in \mathbb{N}$. Then, the generalized $q$-exponential function with respect to $P$ given by
\begin{equation*}
    E_{q,P}(x) = 1 + \sum_{n=1}^{\infty} \frac{x^n}{\prod_{t=1}^n P(q^t)}.
\end{equation*}
If $P(X) = X-1$, $E_{q,P}(x) = E_q(x)$, the $q$-exponential function. Note that $E_{q,P}(x)$ is a basic hypergeometric series, as defined in \cite{gasper-rahman}. \\


With these notations, the first result of this paper is stated below. 
\begin{theorem}\label{thm:most-general-poly}
Let $q$ be such that $\pm q$ is a PV number. Let $P(X)= L_{D} X^{D}+\cdots+ c_{d} X^{d}\in\mathbb{Z}[X]$ be a non-constant polynomial with $P(q^t) \neq 0$ for $t \geq 1$, $d \leq D$ and $L_D \, c_d \neq 0$. Let $\alpha_1,\ldots,\alpha_m$ be non-zero algebraic integers in $\Q(q)$ satisfying  
\begin{equation}\label{eqn:cdn-for-thm}
|c_d|^{n_q-1} \, \max\left\{|\alpha_1|, \, |\alpha_2|, \, \cdots, \, |\alpha_m| \right\} \, \prod_{l=2}^{n_q} \max \left\{ 1, \, \left| \sigma_l(\alpha_1) \right|, \, \left| \sigma_l(\alpha_2) \right|, \, \cdots, \left| \sigma_l(\alpha_m) \right| \right\} < |q|^{D}.
\end{equation}
Suppose that $\alpha_{k_1}/\alpha_{k_2}$ is not a root of unity for $1 \leq k_1,\, k_2 \leq m$ and $k_1 \neq k_2$. Then the numbers in the set
\begin{equation*}
    S := \left\{ E_{q, P}^{(j)} \left( \alpha_k \right)  \, : \, 1 \leq k \leq m, \, 0 \leq j \leq M \right\} \cup \{ \, 1 \, \}
\end{equation*}
are linearly independent over the field $\Q(q)$.
\end{theorem}

The following is an immediate corollary of this theorem.
\begin{corollary}\label{cor:q-exp}
Let $q$ be such that $\pm q$ is a PV number. Let $\alpha_1,\ldots,\alpha_m$ be non-zero algebraic integers in $\Q(q)$ satisfying  
\begin{equation*}
\max\left\{|\alpha_1|, \, |\alpha_2|, \, \cdots, \, |\alpha_m| \right\} \, \, \prod_{l=2}^{n_q} \max \left\{ 1, \, \left| \sigma_l(\alpha_1) \right|, \, \left| \sigma_l(\alpha_2) \right|, \, \cdots, \left| \sigma_l(\alpha_m) \right| \right\} < |q|.
\end{equation*}
Suppose that $\alpha_{k_1}/\alpha_{k_2}$ is not a root of unity for $1 \leq k_1,\, k_2 \leq m$ and $k_1 \neq k_2$. Then the numbers in the set
\begin{equation*}
    S := \left\{ E_q^{(j)}(\alpha_k)  \, : \, 1 \leq k \leq m, \, 0 \leq j \leq M \right\} \cup \{ \, 1 \, \}
\end{equation*}
are linearly independent over the field $\Q(q)$.
\end{corollary}

In particular, this implies the following about the special functions discussed earlier.
\begin{corollary}\label{cor:irr-exp}
Let $q$ be such that $\pm q$ is a PV number and $\alpha \in \mathcal{O}_{q}$ satisfies 
\begin{equation*}
    0 < \bigg( \min\left\{1, \, | \alpha| \right\} \bigg) \, H_q \left( \alpha \right) < |q|.
\end{equation*}
Then the values $E_q(\alpha), \, L_q(\alpha) \not \in \Q(q)$. In particular, $\zeta_q(1)$ is irrational when $\pm q$ is a PV number. 
\end{corollary}
The irrationality and linear independence of the values of the $q$-logarithm function have been extensively studied by various authors. We refer to \cite{zudilin} for a comprehensive history of the problem and investigation of the values of a generalization of the $q$-logarithm function. The irrationality of $\zeta_q(1)$ when $q$ is an integer was first obtained by Erd\H{o}s \cite{erdos}. More recently, Y. Tachiya \cite[Theorem 2]{tachiya-1} proved that $\zeta_q(1) \not\in \Q(q)$ when $q$ is a PV number, which is also implied by Corollary \ref{cor:irr-exp} above.\\

A special function that is closely related to the study of the $q$-exponential function is the Tschakaloff function, given by
\begin{equation*}
    T_q(x) := 1 + \sum_{n=1}^{\infty} \frac{x^n}{q^{\frac{n(n+1)}{2}}}.
\end{equation*}
In our notation, $T_q(x) = E_{q,I}(x)$, where $I(x) = x$. Thus, Theorem \ref{thm:most-general-poly} implies the following.
\begin{corollary}\label{cor:Tschak-values-cor}
Let $q$ be such that $\pm q$ is a PV number. Suppose that $ \alpha \in \mathcal{O}_{q}$ satisfies
\begin{equation*}
  0 <  \bigg(  {   \min \left\{1, \, |\alpha| \right\} } \bigg) \, {H_q(\alpha)} < |q|.
\end{equation*}
Then the numbers
\begin{equation*}
    1, \, T_q(\alpha), \, T_q^{(1)}(\alpha), \, \cdots, \, T_q^{(m)}(\alpha)
\end{equation*}
are linearly independent over $\Q(q)$.
\end{corollary}

\medskip

It was brought to our notice by the referee that Theorem \ref{thm:most-general-poly} follows from Corollaries 5.1 and 5.2 in \cite{a-m-v}, which require a much weaker condition on the $\alpha_k$'s than in Theorem \ref{thm:most-general-poly}. The authors, M. Amou, T. Matala-Aho and K. V\"{a}\"{a}n\"{a}nen, prove a general result regarding linear independence of values of solutions to $q$-difference equations in \cite{a-m-v}. As such, the techniques necessary to prove the result in \cite{a-m-v} are involved, whereas the proof of Theorem \ref{thm:most-general-poly} provided in this paper follows from relatively elementary considerations. \\

The statements so far were concerned with the independence of values of a single function and its derivatives at several arguments. We now address the question of independence of different cognate functions at the same argument. First, for any $M \in \mathbb{N}$ and for any $q$ with $|q|>1$, we define an arithmetic progression analogue of the $E_q(x)$ as
\begin{equation*}
    E_{q,M}(x) := 1 + \sum_{n=1}^{\infty} \frac{x^n}{[Mn]_q!},
\end{equation*}
which is an entire function. Clearly $E_{q,1}(x) = E_q(x)$ and  
\begin{equation*}
    E_{q,M}(x^M) = 1 + \sum_{\substack{n=1, \\ n \equiv 0 \bmod M}}^{\infty} \frac{x^n}{[n]_q !}.
\end{equation*}
Note that $E_{q,M}$ is not a basic hypergeometric function. \\

For these special functions, we prove the following theorem.
\begin{theorem}\label{thm:gap-denominator-lin-ind}
Let $q$ be such that $\pm q$ is a PV number and $a_1<\cdots<a_k$ be distinct positive integers. Let $\alpha \in \mathcal{O}_{q}$ be such that $1 \leq |\alpha|$ and 
\begin{equation}\label{eqn:cdn-thm-2}
 H_q(\alpha) < |q|^{a_1}.
\end{equation}
Then the numbers 
\begin{equation}\label{eqn:thm-4}
 1,  E_{q,a_1}(\alpha),\ldots, E_{q, a_k}(\alpha)
\end{equation}
are linearly independent over the field $\Q(q)$.
\end{theorem}

The approach in this paper is an adaptation of the proof of Theorem 1.1 in \cite{murty-fib}, which is a slight modification of the argument by Duverney \cite{duverney-value-1}. In essence, it is similar to Fourier's proof of irrationality of the number $e$. The proof of Theorem \ref{thm:most-general-poly} relies on a Diophantine lemma, which is a consequence of the Skolem-Mahler-Lech theorem. The proof of Theorem \ref{thm:gap-denominator-lin-ind} is completed using a recursive elimination argument.\\

\medskip

\section{\bf Proof of the theorems}
\bigskip

An important ingredient in the proofs is the following particular case of the Skolem-Mahler-Lech theorem \cite[Theorem 4.3, page 124]{cor-zan}.
\begin{theorem}\label{skm-thm}
Let $\alpha_1,\ldots,\alpha_k$  be non-zero algebraic numbers such that $\alpha_i/\alpha_j$ is not a root of unity for $1\leq i< j\leq k$. Let $P_1(x),\ldots, P_k(x)$ be non-zero polynomials with algebraic coefficients. Then, there are only finitely many natural numbers $n$ satisfying 
\begin{equation*}
P_1(n) \, \alpha^n_1 + \cdots + P_k(n) \, \alpha^n_k=0.
\end{equation*}
\end{theorem}
This is immediate from the Skolem-Mahler-Lech theorem since the sequence $P_1(n) \, \alpha^n_1 + \cdots + P_k(n) \alpha^n_k$ is a non-degenerate recurrence sequence if $\alpha_i/\alpha_j$ is not a root of unity.
\bigskip

\subsection{Proof of Theorem \ref{thm:most-general-poly}}
To begin with, let $f_j(x) := x^j E_{q,P}^{(j)}(x)$ for each $0 \leq j \leq M$. Observe that the result follows if we show that $1$ and the values $f_j(\alpha_k)$ are $\Q(q)$-linearly independent for $0 \leq j \leq M$ and $1 \leq k \leq m$. Indeed, suppose that $\xi_0$ and $\xi_{j,k}$ are algebraic numbers in $\Q(q)$ for $1 \leq k \leq m$ and $0 \leq j \leq M$, not all zero, such that
\begin{equation*}
   \xi_0 + \sum_{j=0}^M \sum_{k=1}^m \xi_{j,k} \, E_{q,P}^{(j)}(\alpha_k)  = 0.
\end{equation*}
Then we obtain the non-trivial linear relation
\begin{equation*}
 \xi_0 + \sum_{j=0}^M \sum_{k=1}^m \frac{\xi_{j,k}}{\alpha_k^j} \, f_j(\alpha_k) = 0,
\end{equation*}
which again has coefficients in $\Q(q)$. Thus it suffices to establish the linear independence of $f_j(\alpha_k)$'s over $\mathbb{Q}(q)$.\\

\medskip

Let $r_0(X) = 1$. For $1 \leq j \leq M$, let $r_j(X) := X (X-1) \cdots (X - j + 1)$. Then 
\begin{equation*}
    f_j(x) = \sum_{n=j}^{\infty} \frac{r_j(n) \, x^n}{ \prod_{t=1}^n P(q^t)} = \sum_{n=1}^{\infty} \frac{r_j(n) \, x^n}{ \prod_{t=1}^n P(q^t)},
\end{equation*}
as $r_j(n) = 0$ for $0 \leq n \leq j-1$. Now suppose $\lambda_0$ and $\lambda_{j,k} \in \Q(q)$ are such that
\begin{equation*}
    \lambda_0 + \sum_{j=0}^M \sum_{k=1}^m \lambda_{j,k} \, f_j(\alpha_k) = 0.
\end{equation*}
Without loss of generality, we can assume that $\lambda_0$ and $\lambda_{j,k}$'s are algebraic integers. For $1 \leq k \leq m$ let $A_k(X) := \sum_{j=0}^M \lambda_{j,k} \, r_j(X)$. Then using the definition of $E_{q,P}(x)$, we get
\begin{align*}
    \widetilde{\lambda_0} + \sum_{n=1}^{\infty} \frac{\sum_{k=1}^m A_k(n) \, \alpha_k^n}{\prod_{t=1}^n P(q^t)} = 0,
\end{align*}
where $\widetilde{\lambda_0} = \lambda_0 + \sum_{k=1}^m \lambda_{0,k}$.

\medskip

Let $N$ be a sufficiently large positive integer. We
truncate the above infinite sum at $N$ and clear denominators to obtain
\begin{align}\label{eqn:X-N-truncate-poly}
  \widetilde{\lambda_0}  \left( \prod_{t=1}^{N} P(q^t) \right)  + \sum_{n=1}^N \left(\sum_{k=1}^m A_k(n)\,\alpha_k^n \  \right)\,  \prod_{t=n+1}^{N} P(q^t)  = - \left( \prod_{t=1}^N P(q^t) \right) \sum_{n=N+1}^{\infty} \frac{\sum_{k=1}^m A_k(n)\,\alpha_k^n }{\prod_{t=1}^n P(q^t)}.
\end{align}
Denote the left hand side of the above equation as $X_N$, namely,
\begin{equation}\label{eqn:X-N-def-poly}
    X_N := \widetilde{\lambda_0}  \left( \prod_{t=1}^{N} P(q^t) \right)  + \sum_{n=1}^N \left(\sum_{k=1}^m A_k(n)\,\alpha_k^n \  \right)\,  \prod_{t=n+1}^{N} P(q^t).
\end{equation}
Then $X_N \in \mathcal{O}_{q}$. Moreover, the right hand side of \eqref{eqn:X-N-truncate-poly} can be written as
\begin{align}
   \bigg( \prod_{t=1}^N P(q^t) & \bigg)  \sum_{n=N+1}^{\infty} \frac{\sum_{k=1}^m A_k(n)\,\alpha_k^n }{\prod_{t=1}^n P(q^t)} \nonumber \\
   & = \frac{\sum_{k=1}^m A_k(N+1) \alpha_k^{N+1}}{P(q^{N+1})} +  \frac{1}{P(q^{N+1})}\mathlarger{\mathlarger{\sum}}_{n=2}^{\infty} \frac{\sum_{k=1}^m A_k(N+n)\,\alpha_k^{N+n} }{\prod_{t=N+2}^{N+n} P(q^t)}. \label{eqn:aux}
\end{align}
For simplicity of notation, let
\begin{equation*}
   \boldsymbol{\alpha} :=  \max \left\{|\alpha_1|, |\alpha_2|, \cdots, |\alpha_m| \right\}
\end{equation*}
Using triangle inequality and the fact that each $A_k(X)$ is a polynomial of degree $M$, we get that for all $\nu>0$,
\begin{align*}
    \left|\sum_{k=1}^m A_k(\nu) \, \alpha_k^{\nu} \right| \leq {\boldsymbol{\alpha}}^{\nu} \, \sum_{k=1}^m |A_k(\nu)| \ll {\nu}^M  {\boldsymbol{\alpha}}^{\nu}.
\end{align*}
Also since $|P(q^t)| \sim |q|^{tD}$ for $t$ sufficiently large, the second term in \eqref{eqn:aux} can be seen to be bounded by
\begin{align*}
    \left| \frac{1}{P(q^{N+1})}\mathlarger{\mathlarger{\sum}}_{n=2}^{\infty} \frac{\sum_{k=1}^m A_k(N+n)\,\alpha_k^{N+n} }{\prod_{t=N+2}^{N+n} P(q^t)} \right|     \ll \frac{{\boldsymbol{\alpha}}^{N+1}}{|P(q^{N+1})|} \mathlarger{\mathlarger{\sum}}_{n=2}^{\infty} {{(n+N)}^M}  \, \cdot \, {\left(\frac{\boldsymbol{\alpha}}{|q|^{DN}} \right)}^{n-1} \, \cdot \, {|q|}^{\frac{- \, D \left( n^2+n-2 \right)}{2}}.
\end{align*}
The above infinite series converges absolutely as $|q|>1$ and the terms decay exponentially. Applying these bounds to the expression in \eqref{eqn:aux} gives
\begin{equation}\label{eqn:first-estimate}
    | X_N | \ll {\frac{\boldsymbol{\alpha}^{N+1}}{|P(q^{N+1})|} }\,N^M,
\end{equation}
where the implied constant depends on $q$, $\alpha_k$'s and the coefficients $\lambda_{j,k}$'s. \\

\medskip 

We now estimate the size of conjugates of $X_N$. Since $\pm q$ is a PV number, $|\sigma_l(q)| < 1$ for $2 \leq l \leq n_q$. From the expression for $X_N$ in \eqref{eqn:X-N-def-poly}, we have for all $n \geq 0$,
\begin{equation*}
    \sigma_l(X_N) = \sigma_l(\widetilde{\lambda_0}) \left( \prod_{t=1}^{N} P \left( {\sigma_l(q)}^t \right) \right)  + \sum_{n=1}^N \left(\sum_{k=1}^m \sigma_l(A_k(n))\,\sigma_l{(\alpha_k)}^n \right)\, \prod_{t=n+1}^{N} P \left( {\sigma_l(q)}^t \right).
\end{equation*}
Observe that 
\begin{equation*}
\left|\prod_{t=n+1}^N P \left( \sigma_l(q^t) \right)\right|= \left|c_d\left(\prod_{t=n+1}^N \sigma_l(q^t)\right)^{d}\right|^{N-n} \, \, \prod_{t=n+1}^N \left| 1+\cdots+\frac{L_D}{c_d} \left( \sigma_l(q^t) \right)^{D-d} \right|.
\end{equation*}
Since $|\sigma_l(q)| < 1$ for all $2 \leq l \leq n_q$, the series
\begin{equation*}
    \sum_{t=1}^{\infty} {\left( \sigma_l \left( q^t \right) \right)}^s
\end{equation*}
is absolutely convergent for all $1 \leq s \leq D-d$. Thus, the infinite product 
$$
\prod_{t=1}^\infty \left| 1+\cdots+\frac{L_D}{c_d} \left( \sigma_l(q^t) \right)^{D-d} \right|
$$
is convergent and we have 
$$
\left|\prod_{t=n+1}^N P \left( {\sigma_l(q)}^t \right) \right| \ll |c_d|^{N-n} \,\, \prod_{t =n+1}^N \left| \left( \sigma_l(q^t) \right)^{d(N-n)} \right| \ll |c_d|^{N-n},
$$
as $|\sigma_l(q)|<1$ for all $2\leq l \leq n_q$. By these observations, we get that 
$$
| \sigma_l(X_N)| \ll|c_d|^{N} \left( 1 \, +  \, \sum_{n=1}^N |c_d|^{-n}\left(\sum_{k=1}^m |\sigma_l(A_k(n))|\,|\sigma_l{(\alpha_k)}|^n \right) \right).
$$
Note that $c_d \in \mathbb{Z}$, Hence, $|c_d| \geq 1$. Now, $\sigma_l \left( A_k(n) \right) = \sum_{j=0}^M \sigma_l \left( \lambda_{j,k} \right) r_j(n)$, which is again a polynomial of degree $M$ in $n$. Putting these bounds together, we deduce that
\begin{equation}\label{eqn:conjugate-estimate-poly}
    \left| \sigma_l(X_N) \right| \ll N^{M+2} \, |c_d|^{N} \, {\bigg(\max \left\{ 1, \, \left| \sigma_l(\alpha_1) \right|, \, \cdots, \left| \sigma_l(\alpha_m) \right| \right\} \bigg)}^N.
\end{equation}
As before, the implied constant above only depends on $q$, $\alpha_k$'s and the $\lambda_{j,k}$'s. \\

\medskip

Multiplying the absolute values of all the conjugates of $X_N$ and the corresponding bounds in \eqref{eqn:first-estimate} and \eqref{eqn:conjugate-estimate-poly}, we derive that
\begin{align*}
    & \hspace{2cm} \prod_{l=1}^{n_q} \left| \sigma_l(X_N) \right| \\
    &\ll \frac{N^{n_q(M+2)-2} \, |c_d|^{(n_q-1)N}\, {\boldsymbol{\alpha}}^{N}}{|P(q^{N+1})|} \,{\left( \prod_{l=2}^{n_q} \max \left\{ 1, \, \left| \sigma(\alpha_1) \right|, \, \cdots, \left| \sigma(\alpha_m) \right| \right\} \right)}^{N}\\
    &\ll N^{n_q(M+2)-2} \,\left(\frac{{\boldsymbol{\alpha}} \, |c_d|^{(n_q-1)}\,{\prod\limits_{l=2}^{n_q} \max \left\{ 1, \, \left| \sigma(\alpha_1) \right|, \, \cdots, \left| \sigma(\alpha_m) \right| \right\}}}{|q|^{D}}\right)^N 
\end{align*}
By the hypothesis \eqref{eqn:cdn-for-thm}, the above bound tends to zero as $N \rightarrow \infty$. In particular, this implies that
\begin{equation*}
   \left| \prod_{l=1}^{n_q}  \sigma_l(X_N) \right| < 1
\end{equation*}
for all $N$ sufficiently large. However, the left hand side above is a power of the norm of an algebraic integer. Note here that $\mathbb{Q}(X_N)$ may be a strict sub-field of $\mathbb{Q}(q)$. Thus, $ \prod_{l=1}^{n_q}  \sigma_l(X_N)$ must be a rational integer for all $N>0$. This is only possible if $X_N = 0$ for \textit{all} $N$ sufficiently large. \\

Therefore, there exists a natural number $N_0$ such that for all $N \geq N_0$, 
\begin{equation*}
    \frac{X_N}{\prod_{t=1}^N P(q^t)} = \widetilde{\lambda_0} + \sum_{n=1}^N \sum_{k=1}^m A_k(n) \, \alpha_k^n = 0.
\end{equation*}
Thus, considering the expression, 
\begin{equation*}
     \frac{X_{N+1}}{\prod_{t=1}^{N+1} P(q^t)} - \frac{X_N}{\prod_{t=1}^N P(q^t)},
\end{equation*}
which equals zero for $N > N_0$, we obtain that
\begin{equation*}\label{eqn:recurrence}
  A_1(N) \, \alpha_1^N + \cdots + A_m(N) \,\alpha^N_m = 0,
\end{equation*}
for all $N > N_0$.  As $\alpha_{k_1}/\alpha_{k_2}$ is not a root of unity, Theorem \ref{skm-thm} applied to the above relation immediately gives that $A_k(N) = 0$ for all $1 \leq k \leq m$ and $N > N_0$. Thus, the polynomials $A_k(X)$ are identically zero. Recall that
\begin{equation*}
    A_k(X) = \sum_{j=0}^M \lambda_{j,k} \, r_j(X), 
\end{equation*}
and $\deg r_j(X) = j$. Since $r_j(X)$ have distinct degrees, $A_k(X)$ is identically zero if and only if $\lambda_{j,k} = 0$ for all $0 \leq j \leq M$ and $1 \leq k \leq m$. This completes the proof of the theorem.
\qed

\bigskip

\subsection{Proof of Theorem \ref{thm:gap-denominator-lin-ind}}
We begin along the same lines as the argument in the proof of Theorem \ref{thm:most-general-poly}.\\

Suppose that the numbers in \eqref{eqn:thm-4} are linearly dependent over $\Q(q)$. Then, there exist algebraic integers $\lambda_0, \lambda_1, \ldots, \lambda_k \in \mathcal{O}_q$ not all zero such that
\begin{equation*}
 \lambda_0+\lambda_1 E_{q,a_1}(\alpha)+\cdots+ \lambda_k E_{q,a_k}(\alpha)=0. 
\end{equation*}
Without loss of generality, we can assume that $\lambda_1 \neq 0$. For otherwise, we can change notation to replace $a_j$ by $a_1$ for the smallest $j \leq k$ for which $\lambda_j \neq 0$ and work out the argument below. \\

Using the definition of the $q$-exponential function, we get
\begin{align}\label{eqn:thm-2-reln}
    \widetilde{\lambda_0} + \sum_{n=1}^{\infty} \frac{\lambda_1\,\alpha^n}{[a_1 n]_{q !}} +\cdots+ \sum_{n=1}^{\infty} \frac{\lambda_k\,\alpha^n}{[a_k n]_{q}!}= 0,
\end{align}
where $\widetilde{\lambda_0} = \lambda_0 + \lambda_1 + \cdots + \lambda_k$. Set $d=\mbox{lcm}\{a_1,\ldots,a_k\}$ and $d_i=d/a_i$. Choose a large positive integer $N$ which is our parameter and  set  $N_i=N d_i$ for $i=1, \, 2,\ldots,k$. With these choices of $N_i$, we have
\begin{equation*}
    a_1 N_1 = a_2 N_2 = \cdots = a_k N_k = d N.
\end{equation*}

Furthermore, for all $i=1, \,2, \, 3,\ldots,k$, 
\begin{align}
 \frac{[dN]_{q}! }{[a_i(N_i+1)]_{q}!} = \frac{[a_i N_i]_{q}! }{[a_i(N_i+1)]_{q}!}&=\frac{(q^{a_i N_i}-1)\cdots(q-1)}{(q^{a_i N_i+a_i}-1)\cdots(q-1)}=\frac{1}{(q^{N d+a_i}-1)\cdots(q^{N d+1}-1)}.\label{eqn:relation2}
\end{align}
\smallskip

\noindent
Now truncate the infinite sums in \eqref{eqn:thm-2-reln} at $N_i$'s and multiply by $[dN]_{q}!$ to get 
\begin{align}
[dN]_{q}!\bigg(\widetilde{\lambda_0} + \sum_{n=1}^{N_1} \frac{\lambda_1\,\alpha^n}{[a_1 n]_{q} !}&  +\cdots+ \sum_{n=1}^{N_k} \frac{\lambda_k\,\alpha^n}{[a_k n]_{q} !} \bigg)  \nonumber \\
& =-[dN]_{q}!\left(\sum_{n=N_1+1}^{\infty} \frac{\lambda_1\,\alpha^n}{[a_1 n]_{q} !} +\cdots+ \sum_{n=N_k+1}^{\infty} \frac{\lambda_k\,\alpha^n}{[a_k n]_{q} !}\right). \label{eqn:truncate-thm3}
\end{align}
Let $X_N$ denote the left hand side of the above equation, i.e.,
\begin{equation}\label{eqn:X-N-1}
X_N:= [dN]_{q}!\left(\widetilde{\lambda_0} + \sum_{n=1}^{N_1} \frac{\lambda_1\,\alpha^n}{[a_1 n]_{q} !} +\cdots+ \sum_{n=1}^{N_k} \frac{\lambda_k\,\alpha^n}{[a_k n]_{q} !}\right).    
\end{equation}
Since $[dN]_q! = [a_iN_i]_q!$ for $1 \leq i \leq k$, $X_N$ is an  algebraic integer in $\mathcal{O}_q$. We now estimate the right hand side of \eqref{eqn:truncate-thm3}. Note that by an argument similar to the one in Theorem \ref{thm:most-general-poly} and using \eqref{eqn:relation2}, we can deduce
\begin{align*}
    \left|[dN]_q! \sum_{n = N_j+1}^{\infty} \frac{\alpha^n}{[a_jn]_q!} \right| \ll \frac{{|\alpha|}^{N_j}}{|q|^{a_jdN}} \ll {\left( \left| \frac{\alpha^{d_j}}{q^{a_jd}} \right| \right)}^N,
\end{align*}
since $N_j = Nd_j$. Since $a_1 < a_2 < \cdots < a_k$, $d_1 > d_2 > \cdots > d_k$ and $|\alpha| \geq 1$, we derive that
\begin{equation}\label{eqn:estimate-X-N-1}
    |X_N| \ll {\left(  \left| \frac{\alpha^{d_1}}{q^{a_1d}} \right| \right)}^N.
\end{equation}
By the same argument as in the proof of Theorem \ref{thm:most-general-poly}, we obtain the following estimate for the conjugates of $X_{N}$: 
\begin{equation}\label{eqn:estimate-conjugate-X-N-1}
 |\sigma_l(X_N)|\ll  N_1  \, {\left( \max\{1, |\sigma_l(\alpha)|\} \right)}^{N_1}. 
\end{equation}
As before, the implied constant above only depends on $q$, $a_i$'s and the $\lambda_{j,k}$'s. Multiplying the bounds \eqref{eqn:estimate-X-N-1} and \eqref{eqn:estimate-conjugate-X-N-1} on the absolute values of all the conjugates of $X_N$ and noting that $|\alpha| \geq 1 $, we derive that
\begin{align*}
    \prod_{l=1}^{n_q} \left| \sigma_l(X_N) \right| \ll N_1^{n_q-1} \,  \, { \left( \frac{|\alpha| \prod\limits_{l=2}^{n_q} \max \left\{ 1, \, \left| \sigma_l(\alpha) \right| \right\} }{|q|^{a_1^2}} \right) }^{d_1N} .
\end{align*}
By \eqref{eqn:cdn-thm-2}, the right hand side above tends to zero as $N \rightarrow \infty$. However, the left hand side is a rational integer since it is a power of the norm of an algebraic integer. Therefore, there exists a natural number $N_0$ such that for \textit{all} $N> N_0$, $X_N = 0$, which in turn implies that $X_N=X_{N+1}=0$. Consequently,  
\begin{equation*}
  \widetilde{\lambda_0} + \sum_{n=1}^{N d_1} \frac{\lambda_1\,\alpha^n}{[a_1 n]_{q}!} +\cdots+ \sum_{n=1}^{N d_k} \frac{\lambda_k\,\alpha^n}{[a_k n]_{q}!}=0  
\end{equation*}
and 
\begin{equation*}
  \widetilde{\lambda_0} + \sum_{n=1}^{N d_1+d_1} \frac{\lambda_1\,\alpha^n}{[a_1 n]_{q}!} +\cdots+ \sum_{n=1}^{N d_k+d_k} \frac{\lambda_k\,\alpha^n}{[a_k n]_{q}!}=0  
\end{equation*}
for all $N > N_0$. By subtracting these two relations, we get that 
\begin{equation}\label{eqn:X-N-d}
 \lambda_1 \sum_{n=N d_1+1}^{N d_1+d_1} \frac{\alpha^n}{[a_1 n]_{q}!} +\cdots+ \lambda_k \sum_{n=N d_k+1}^{N d_k+d_k} \frac{\alpha^n}{[a_k n]_{q}!}=0 
\end{equation}
for all $N > N_0$. Note that for $1 \leq j \leq k$, we have $Nd+a_j \leq a_j n \leq Nd + d$ in the above sums. Therefore, the term 
\begin{equation*}
    \frac{\alpha^{Nd_k+1}}{[Nd+a_1]_q!}
\end{equation*}
divides each term in the above relation. Factoring this out, we obtain
\begin{align}
    & \lambda_1 \left(\alpha^{N(d_1 - d_k)} + \frac{\alpha^{N(d_1 - d_k)+1}}{(q^{Nd+2a_1}-1)\cdots(q^{Nd+a_1 + 1}-1)} + \cdots + \frac{\alpha^{N(d_1 - d_k)+d_1 -1}}{(q^{Nd+d}-1)\cdots(q^{Nd+a_1 + 1}-1)} \right) + \nonumber \\
    &  \lambda_2 \left(\frac{\alpha^{N(d_2 - d_k)}}{(q^{Nd+a_2}-1)\cdots(q^{Nd+a_1 + 1}-1)}  + \cdots + \frac{\alpha^{N(d_2 - d_k)+d_2 -1}}{(q^{Nd+d}-1)\cdots(q^{Nd+a_1 + 1}-1)}\right) \nonumber \\
    & + \nonumber \\
    & \vdots \nonumber \\
    & + \nonumber \\
    & \lambda_k \left( \frac{1}{(q^{Nd+a_k}-1)\cdots(q^{Nd+a_1 + 1}-1)}  + \cdots + \frac{\alpha^{d_k-1}}{(q^{Nd+d}-1)\cdots(q^{Nd+a_1 + 1}-1)}\right) = 0. \label{eqn:rel-a-1}
\end{align}\\
Now, for $1 \leq j \leq k$ and $0 \leq l \leq d_j-1$, the absolute value of a general term is of the form
\begin{equation*}
    \left| \frac{\alpha^{N(d_j - d_k) + l}}{(q^{Nd+(l+1)a_j}-1)\cdots(q^{Nd+a_1+1}-1)}\right| \ll {\left| \frac{\alpha^{d_1 - d_k}}{q^{\delta d}} \right|}^N
\end{equation*}
except for $j=1$ and $l=0$, with $\delta = \min \{a_1, \, a_2 - a_1\}$. Since $1 \leq \delta$, this implies that each term in the above relation is $\ll {|\alpha^{d_1 - d_k}/q^d|}^N$. By \eqref{eqn:cdn-thm-2}, this quotient is less than $1$ as $1 \leq |\alpha| < |q|^{a_1}$. Hence, taking the limit as $N \rightarrow \infty$, all terms in \eqref{eqn:rel-a-1} tend to zero except the very first one, that is, $\alpha^{N(d_1-d_k)}$. This implies that $\lambda_1 = 0$, which is a contradiction. This completes the proof.\\

\qed

\bigskip

\section*{\bf Acknowledgments}
This work was initiated when the second author was visiting Chennai Mathematical Institute during July 2022. He extends heartfelt gratitude to the institute for its' hospitality and support. We also thank Prof. Ram Murty for helpful suggestions to improve the exposition of the paper. We are especially grateful to the referee for a careful reading and insightful comments about the earlier draft of this paper, and thank them for making us aware of the results in \cite{a-m-v}.

\end{document}